\documentclass[12pt]{article}
\input{epsf}
\usepackage{amsmath}
\usepackage{amsfonts}
\usepackage{amscd}
\newtheorem{Def}{Definition}

\newtheorem{Thm}[Def]{Theorem}

\def\block{\hbox{${\vcenter{\vbox{\hrule height 0.4pt\hbox{\vrule 
width 0.4pt 
height 6pt \kern 5pt\vrule width 0.4pt}\hrule height 0.4pt }}}$}}

\def\Hom{\mbox{Hom}}

\def\H{{\cal H}}

\begin{document}

\title{2-Categories, 4d state-sum models and gerbes\footnote{Talk presented 
by the first author at III Encontro F\'{o}rum Internacional de 
Investigadores Portugueses, 
Faro, April 7-10, 2001.}}
\author{$\mbox{Marco Mackaay}^{1,3}$
\and $\mbox{Roger Picken}^{2,3}$}
\date{April 30, 2001}
\maketitle

\noindent In this article we focus on the third word in the title of our talk and 
on our motivation for getting involved with it. 

The simplest state-sum model is the Dijkgraaf-Witten (DW) model~\cite{DW90}. 
In its most elementary form the DW-model associates to a smooth closed 
connected finite-dimensional manifold $M$ the following number: 
$$\#\Hom\left\{\pi_1(M),G\right\},$$ where $G$ is a finite group. If we want 
to understand the differential geometry behind the DW-model we have to give 
up the finiteness of $G$ of course. If $G$ is a Lie group, we can ask 
ourselves what geometric objects correspond to smooth homomorphisms from 
$\pi_1(M)$ to $G$ (we will not explain here what we mean by smoothness exactly). The 
answer is well known: {\em principal $G$-bundles with flat connections}. We 
explain this in some more detail. 
Let $\left\{U_i\right\}$ be a covering of $M$ by open sets such that all 
intersections 
$U_{i_1\ldots i_p}=U_{i_1}\cap\cdots\cap U_{i_p}$ are contractible. We 
present a principal $G$-bundle, $P$, by its transition functions 
$g_{ij}\colon U_{ij}\to G$, which satisfy $g_{ji}=g_{ij}^{-1}$ and the 
cocycle condition 
$$g_{ij}g_{jk}g_{ik}^{-1}=1\quad\mbox{on}\ U_{ijk}.$$
A connection, $\cal A$, in $P$ can be defined in terms of local 1-forms,  
$A_i$ on $U_i$, with values in the Lie algebra of $G$, which satisfy
$$A_j-g_{ij}^{-1}A_ig_{ij}=g_{ij}^{-1}dg_{ij}\quad\mbox{on}\ U_{ij}.$$ 
Given a loop $\ell$ in 
$M$ one can define the {\em holonomy}, $\H(\ell)\in G$, of $\cal A$ around 
$\ell$. In general the holonomies around two homotopic loops are different. 
However, if the two loops are {\em thin homotopic}, then the holonomies are 
equal, as was first remarked by Barrett~\cite{Ba91}. 
There are several ways to define thin homotopy mathematically~\cite{Ba91,
CP94}. 
For the purpose of this article it suffices to give the intuitive idea: 
\begin{Def}{\rm(sketch)} Two loops are thin homotopic if there exists a 
homotopy between them whose image has no area. 
\end{Def}
All homotopies involved in the standard proof that 
$\pi_1(M)$ is a group are thin as a matter of fact, so we can define the 
{\em thin fundamental group} of $M$, denoted $\pi_1^1(M)$, by dividing out 
the set of loops only by thin homotopies. Note that $\pi_1(M)$ is a 
quotient of $\pi_1^1(M)$. Thus $\cal A$ gives rise to 
a {\em holonomy homomorphism} $$\H\colon \pi_1^1(M)\to G,$$
which is smooth in a technical sense~\cite{Ba91,CP94}. Barrett~\cite{Ba91} 
(see~\cite{CP94} for a proof of the analogous statement using a different 
definition of thin homotopy) proved that 
there is a converse statement:
\begin{Thm} Given a smooth homomorphism $\H\colon \pi_1^1(M)\to G$, there is 
a principal $G$-bundle with connection, unique up to equivalence, 
whose holonomy homomorphism is equal to $\H$. 
\end{Thm}
A connection is flat precisely when the corresponding $\H$ factors through 
the ordinary $\pi_1(M)$. 

Let us now assume that $\pi_1(M)=0$. The next state-sum model that we consider 
associates to $M$ the number
$$\#\Hom\left\{\pi_2(M),H\right\},$$
where $H$ is a finite Abelian group. This is a special case of the 
Yetter model~\cite{Ye93}, which involves the homotopy 2-type of $M$. The right 
algebraic framework for the Yetter model is that of 2-categories, which is 
how the first word in the title of our talk enters the picture. For more 
information about 2-categories and 4-dimensional state-sum models see the two 
papers~\cite{Ma99,Ma00} by the first author and references therein. 
The question now arises whether we can 
understand the maps $\pi_2(M)\to H$ in an analogous differential geometric 
way. Let $H=U(1)$. We see immediately 
that the answer cannot be found in the framework of bundles and connections, 
because we need some geometric structure that gives rise to holonomy around 
surfaces rather than loops. It is known 
(well-known would be an over-statement) that {\em gerbes} with 
{\em gerbe-connections} give rise to such holonomies~\cite{Bry93}. 
A gerbe, $\cal G$, can be defined by functions on triple intersections, 
$h_{ijk}\colon U_{ijk}\to U(1)$, which satisfy 
$h_{\sigma(i)\sigma(j)\sigma(k)}=h_{ijk}^{\epsilon(\sigma)}$, for any 
$\sigma\in S_3$, and the next order cocycle condition:
$$h_{ijk}h_{ijl}^{-1}h_{ikl}h_{jkl}^{-1}=1\quad\mbox{on}\ U_{ijkl}.$$
A gerbe-connection, $\cal B$, in $\cal G$ can be defined by $1$-forms, $A_{ij}$ on $U_{ij}$, and 
$2$-forms, $F_i$ on $U_i$, all with values in $i\mathbb{R}$, such that 
$A_{ji}=-A_{ij}$ and 
$$A_{ij}+A_{jk}-A_{ik}=h_{ijk}^{-1}dh_{ijk}\quad\mbox{on}\ U_{ijk},$$
$$F_j-F_i=dA_{ij}\quad\mbox{on}\ U_{ij}.$$
Gerbes were first defined by Giraud~\cite{Gi71}. The standard reference 
nowadays is Brylinski's book~\cite{Bry93}. The properties of gerbes and gerbe-connections are 
analogous to 
those of line-bundles; e.g. the curvature of a gerbe-connection, 
$\Omega\vert_{U_i}=dF_i$, is a closed 
integral 3-form, the cohomology class of which classifies the gerbe 
up to equivalence, and every closed integral 3-form is the curvature of a 
certain gerbe-connection in a certain gerbe. One can define the 
gerbe-holonomy~\cite{Bry93}, $\H(s)$, of $\cal B$ around any smooth map 
$s\colon S^2\to M$, henceforth referred to as a {\em 
$2$-loop}. 
As for ordinary connections one can show~\cite{MP01} that the gerbe-holonomies 
around two thin homotopic $2$-loops are equal. 
\begin{Def}{\rm(sketch)} We say that two $2$-loops are 
thin homotopic if there exists a homotopy between them whose image has no volume.
\end{Def}
The 
higher dimensional thin homotopy groups were first defined by Caetano and 
Picken in \cite{CP98}, where one can find the technical definition. Thus $\cal B$ gives rise to a smooth {\em gerbe-holonomy homomorphism} 
$$\H\colon\pi_2^2(M)\to U(1),$$
where $\pi_2^2(M)$ is the {\em thin second homotopy group}. In~\cite{MP01} we 
proved that there is a converse statement:
\begin{Thm}\label{ge} Assume that $M$ is simply-connected. Given a smooth homomorphism 
$\H\colon\pi_2^2(M)\to U(1)$, 
there exists a gerbe with gerbe-connection, unique up to equivalence, 
whose holonomy map is equal to $\H$.
\end{Thm}
A gerbe-connection is flat precisely when its holonomy map factors through 
the ordinary $\pi_2(M)$, so we have achieved our goal of understanding 
the differential geometry of the Yetter model. The proof 
of Thm.~\ref{ge} is fairly straightforward. If one does not assume 
that $M$ is simply-connected 
the analogous statement and its proof involve the less familiar 
mathematics of Lie (2-)groupoids. We believe that 
this case is very interesting because it fuses ideas from category theory 
and geometry into something that can best be called 
{\em categorical geometry}. 
\vskip0.2cm
\centerline{\bf Acknowledgements}
This work was supported by {\em Programa Operacional
``Ci\^{e}ncia, Tecnologia, Inova\c{c}\~{a}o''} (POCTI) of the
{\em Funda\c{c}\~{a}o para a Ci\^{e}ncia e a Tecnologia} (FCT),
cofinanced by the European Community fund FEDER.  
The first author is presently on leave from the Universidade do Algarve 
and working with a postdoctoral fellowship from FCT at the University of 
Nottingham (UK).

1. Dep. de Matem\'{a}tica, Univ. do Algarve, Faro, Portugal. {\it Email}: 
mmackaay@ualg.pt.\\
2. Dep. de Matem\'{a}tica, Instituto Superior T\'{e}cnico, Lisboa, 
Portugal. {\it Email}: rpicken@math.ist.utl.pt.\\
3. Centro de Matem\'{a}tica Aplicada, 
Instituto Superior T\'{e}cnico, Lisboa, Portugal. 
\end{document}